# Linguistic Approach to Time Series Forecasting

Dmytro Lande [1, 2], Volodymyr Yuzefovych [1] and Yevheniia Tsybulska [1]

[1] *Institute for Information Recording of National Academy of Sciences of Ukraine, Kyiv, Ukraine*
[2] *National Technical University of Ukraine "Igor Sikorsky Kyiv Polytechnic Institute", Kyiv, Ukraine*

**Abstract**
This paper proposes methods of predicting dynamic time series (including non-stationary ones) based on a linguistic approach, namely, the study of occurrences and repetition of so-called N-grams. This approach is used in computational linguistics to create statistical translators, detect plagiarism and duplicate documents. However, the scope of application can be extended beyond linguistics by taking into account the correlations of sequences of stable word combinations, as well as trends. The proposed methods do not require a preliminary study and determination of the characteristics of time series or complex tuning of the input parameters of the forecasting model. They allow, with a high level of automation, to carry out short-term and medium-term forecasts of time series, characterized by trends and cyclicality, in particular, series of publication dynamics in content monitoring systems. Also, the proposed methods can be used to predict the values of the parameters of a large complex system with the aim of monitoring its state, when the number of such parameters is significant, and therefore a high level of automation of the forecasting process is desirable. A significant advantage of the approach is the absence of requirements for time series stationarity and a small number of tuning parameters. Further research may focus on the study of various criteria for the similarity of time series fragments, the use of nonlinear similarity criteria, the search for ways to automatically determine the rational step of quantization of the time series.

## 1. Introduction

Time series forecasting is a relevant and required technique in various application areas. It should be noted that there are no universal solutions, as evidenced by the presence more than 200 methods of solving this problem, according to some estimates.

There are many forecasting tasks where the forecast should be obtained quickly and with minimal involvement of analysts to organize and implement this process. For example, this approach is applied in solving the problem of parallel predicting of information flow intensity on various topics [1] or in predicting the values of numerous parameters of some large complex system to monitor its condition [2]. In such cases, in addition to the requirement of forecast accuracy, there is also a requirement for a minimum number of the forecasting model parameters to be configured and a simple configuration of those parameters in limited number of iterations. Of particular note are the dynamic time series (where the change in values is mainly caused by the dynamics of the observed process, as opposed to random factors), and in the general case – non-stationary ones (where the mathematical expectation, mean, covariance of some their sub-series are not constant). Recently the method of forecasting time series ARIMA (Autoregressive Integrated Moving Average), which involves the construction of an autoregressive integrated moving average model, has become widespread. In addition, the seasonal modification of the ARIMA model (SARIMA) and a number of similar ones are used, in particular Generalized ESD (Extreme Studentized Deviant), Seasonal Hybrid ESD [3]. Another popular approach is the Sornette wave method [4]. It requires a large number of settings and can be used to analyze non-stationary time series, even economic crisis data [5].

These methods are considered very flexible, as they allow "fine" adjustment of the time series model parameters. However, they require prior detailed study of the series and significant knowledge and experience on the part of the analyst performing the model forecast, and hence additional time and human costs. The solution to this problem is development and use of specialized software products with powerful functionality for automated forecasting. However, a common weakness of such systems is poor ability to explain the resulting forecast. It is difficult for analyst to understand why this and not another forecast was obtained, so in fact, it makes him an outside observer of the process, who must believe obtained result almost without explanation. In addition, such software products are difficult and sometimes impossible to integrate into other automated systems, where the time series forecasting is an intermediate step to achieve some more general goal.

However, there are a number of simpler time series forecasting methods which do not require complex settings and at the same time allow to obtain fairly accurate forecast. Such methods include the so-called exponential methods, which involve, for example, construction of a simple moving average model, Brown's model (simple exponential smoothing) or Holt's model (double exponential smoothing) [6]. As opposed to simple exponential smoothing (Brown's model) or moving average with a given window, Holt's model involves smoothing not only the time series values, but also its trend to make a forecast. This two-parameter model of the time-series $X = x_1, x_2, \ldots, x_k, \ldots, x_K$ of length $K$ includes three equations

$$\begin{cases} \hat{x}_k = (1-\xi)x_k + \xi(\hat{x}_{k-1} + \hat{T}_{k-1}); \\ \hat{T}_k = (1-\varphi)(\hat{x}_k - \hat{x}_{k-1}) + \varphi\hat{T}_{k-1}; \\ \hat{x}_{k+P} = \hat{x}_k + P\hat{T}_k \end{cases} \quad (1)$$

where $\hat{x}_k, \hat{x}_{k-1}$ – the estimates of the series smoothed values at the $k$-th and $k$-1 smoothing steps;

$\xi$ – smoothing coefficient of the series values ($\xi$ varies from 0 to 1);

$T_k, \hat{T}_{k-1}$ – the estimates of the time series values changing trend at the $k$-th and $k$-1 smoothing steps;

$\varphi$ – trend smoothing coefficient ($\varphi \in [0,1]$);

$P$ – the number of future series values to predict (bias time).

The difficulty of using this model is the need to select two smoothing parameters ($\xi$ and $\varphi$), that jointly affect the forecasting result. Note that such selection cannot be carried out consistently, because the equations for the time series values and trend smoothing are related by the same estimates $\hat{X}$. That is, the values of both coefficients affect both equations.

Also there is an advanced method – Holt-Winter's method [7] with three input parameters, that includes additional time series seasonality parameter. This method is used if the series contains cyclic variations. In this case, the time series should be pre-checked for such seasonality, for example, by analyzing its autocorrelation or partial autocorrelation function.

Thus, the application of even simple forecasting methods based on non-stationary time series, require a fairly active role of the analyst, that complicates the process of their automated use to make a quality forecast.

## 2. Purpose of the study

The purpose of this study is to develop methods for predicting dynamic time series, including non-stationary ones, which contain a minimum number of input parameters and at the same time do not require complex settings.

This paper proposes a new, so-called linguistic approach to forecasting, on the basis of which several methods are built. Using this approach requires setting only two parameters that are not directly related, are adjusted separately and do not require in-depth analysis of the time series. The proposed approach allows forecasting of both stationary and non-stationary time series, does not provide the time series analysis for seasonality, but automatically takes it into account when forming the forecast values.

# 3. Linguistic and linguo-correlation methods for the dynamical time series forecasting

Linguistic forecasting method is based on N-gram method, which is used for natural language processing, in particular, to predict the next words in the statement, if all previous ones are known. This task is often solved in computer translators and systems for duplicate and plagiarism detection [8]. It is believed that the conditional probability of the next word depends on the previous words and their sequences. As is known, the assumption that the probability and time series values are determined by previous values is also used to predict the time series. The idea of applying the N-gram method to predict time series is based on a model assumption, that the time series can be considered as some meaningful text, and the new value of this series (the next word) can be determined by analyzing of the previous "text".

Thus, to determine an unknown "phrase" (series prediction) of the length $P$ ($\hat{X}_P$) it would search in the full "text" (time series) $X$ a "phrase" $X_N$ of the length $N$ ($X_N = x_{k+1}, x_{k+2}, …, x_N$), that is identical (similar) to a "phrase" $X_K$, which immediately precedes the forecast moment ($X_K = x_{K-N}, x_{K-N+1}, …, x_K$), and to assume that $P$ "words" that are immediately after the found "phrase" $X_N$ – $X_P = x_{N+1}, x_{N+2}, …, x_{N+P}$, form the desired "phrase" (predicted values of the series), so $\hat{X}_P = X_P$.

Element-by-element difference (the difference criterion) of series values can be used to compare phrases and in the case of their complete coincidence it is equal to 0. However, it is obvious that a full-text search for a "phrase" $X_N$, that exactly matches a "phrase" $X_K$, will not necessarily be successful. Note that increasing the probability of finding identical phrases in the text can be achieved by reducing the number of possible different words in the text (language vocabulary) or by reducing their length. Otherwise, instead of searching for identical phrases, you need to look for similar ones, which can obviously worsen the quality of this forecast. The criterion for the phrases similarity can be the minimal element-by-element difference of the series values, or the maximum correlation coefficient (correlation criterion). It is calculated by expression

$$R_{X_N X_K} = \frac{\sum_{k=1}^{N}(x_k - \bar{x})(y_k - \bar{y})}{\sqrt{\sum_{k=1}^{N}(x_k - \bar{x})^2 \sum_{k=1}^{N}(y_k - \bar{y})^2}}, \quad (2)$$

where  $x_k, y_k$ – values of the series $X_N$ and $X_K$, respectively;
$\bar{*}$ – averaging (mathematical expectation estimation).

If the criterion (2) is used we talk about linguo-correlation forecasting method.

In accordance with the above, an algorithm for solving a forecast problem of the time series $X$ for the bias time $P$ is as follows.

From the time series $X$, $N$ values that immediately precede the moment of forecast are selected successively ($K$). The value $N$ is the first parameter of forecasting model, which must be set by the analyst (researcher) before the forecasting process.

Instead of the parameter $N$ in practice it is convenient to use another parameter $M$ – a real number that specifies the multiplicity of the value $P$. This parameter is selected from product $N = [M \cdot P]$, where $M$ takes the greater value, the smaller the value $P$, and $[*]$ is the rounding function to the nearest larger integer.

As studies have shown, M should be set based on the following conditions.
- If $P \geq 5$, then $1 \leq M \leq 2$, while the larger $P$, the smaller $M$;
- If $P < 5$, then $2 < M \leq 5$, while the smaller $P$, the larger $M$.

Further, to increase the probability of finding identical phrases, the considered time series is subjected to the quantization procedure – the splitting the range of possible series values into a finite number of levels with rounding of actual values to the nearest levels. The purpose of this procedure is to reduce all possible values of the time series to their limited list – "words from the text vocabulary". The quantization step should be maximum, but such that allows to keep the basic changing dynamics of the series values, so the number of different "words" of the language should be the minimum required.

The example of quantization procedure is shown in **Figure 1**: Results of the time series quantization

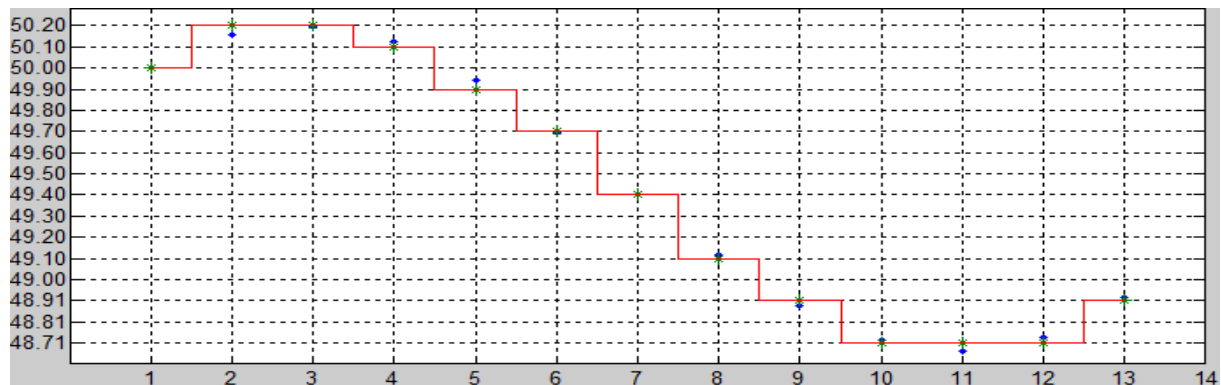

**Figure 1**: Results of the time series quantization

Presented in Figure 1 fragment of the time series as a quantization procedure result is described by nine different "words" arranged in a certain sequence.

In fact, the quantization step or the number of artificial language "words" ($S$), is the second (and final) parameter to be set by the analyst (researcher) to execute forecasting procedure. The peculiarity of this forecasting method is that minor changes in the parameters $N$ and $S$ do not significantly affect the prediction outcome. They are easy to fine tune the final part of the available time series as a test to check the quality of the forecast.

As an example of linguistic method application, consider the forecasting problem at the bias time $P=20$ of some deterministic nonstationary time series of real numbers with the length $K = 100$, and with a maximum range of values 4 units, which is characterized by visually noticeable seasonality, but has no off-season changes in the trend. Such time series can be formed by any periodic function values, such as sinusoid. In accordance with the above rule, set the value $M = 1$ (then $N=P$) and $S=32$ (the quantization step – 4/30=0,125). This time series and forecasting result by linguistic method are shown in the Figure 2.

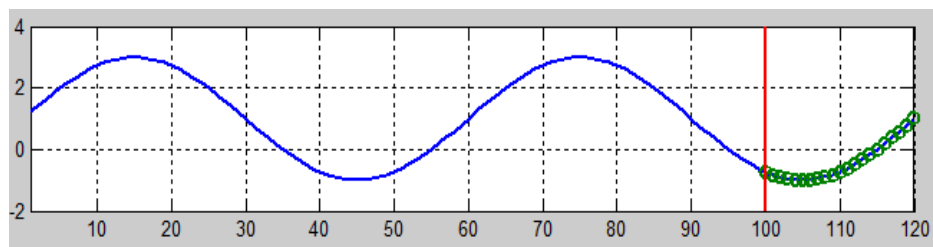

**Figure 2**: Forecasting results of the time series without trend by linguistic method

It is seen that the proposed method in this case allows to solve the forecasting problem quite well, even without calculations of forecast errors.

The use of time series quantization also increases the resistance of the forecast to minor random perturbations of the series values without the use of smoothing procedures. **Figure 3** shows the forecasting results of the same time series, but with random perturbations: uniform distribution law of perturbations in the range from 0.15 to 0.15 units, with zero mathematical expectation.

Consider the more complex problem of predicting a deterministic nonstationary time series of real numbers, which is characterized by seasonality and has a linear trend. Obviously, in this case the complete (infallible) coincidence of the phrases $X_N$ and $X_K$ may not be performed. Therefore, before comparing the phrases $X_N$ and $X_K$, it is necessary to remove their trends. The trend of the time series can be calculated by the least squares method. For a linear trend described by equation $y = Bx + A$ the coefficients $B$ and $A$ are calculated by known expressions

$$\begin{cases} A = \bar{y} - B\bar{x} \\ B = \dfrac{\overline{xy} - \bar{x} \cdot \bar{y}}{\overline{x^2} - (\bar{x})^2} \end{cases} \quad (3)$$

where  $x$ – serial number of the time series value;
  $y$ – time series value.

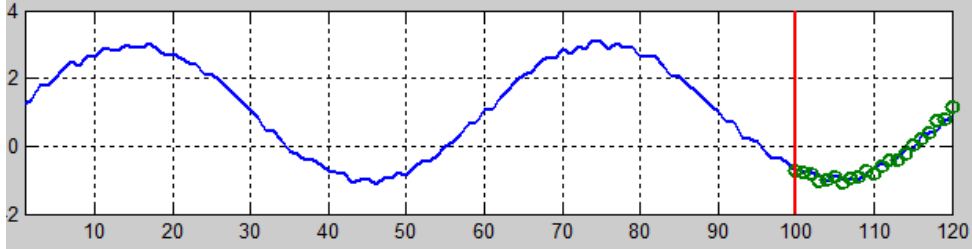

**Figure 3:** Forecasting results of the time series with random perturbations by linguistic method

Next, after finding the "phrases" $X_N$ and $X_P$, to the values of $X_P$ must add a trend that characterizes $X_K$.

**Figure 4** shows forecasting results of deterministic nonstationary time series of real numbers ($K$=100, $P$=20, $M$ = 1 ($N$=$P$) та $S$=30), which is characterized by seasonality and has a "global" linear trend.

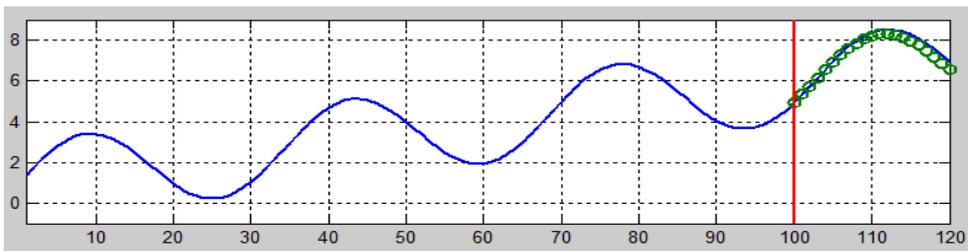
a) Linguistic method applying

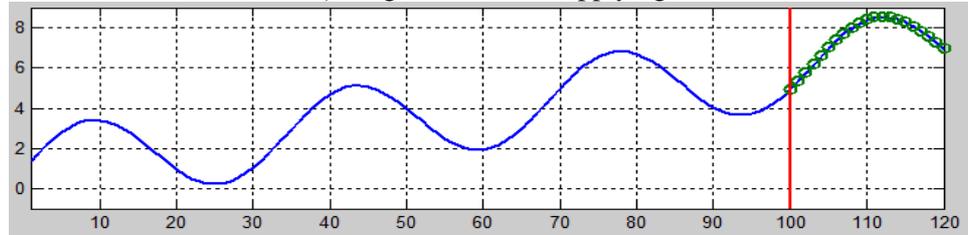
b) Linguo-correlation method applying

**Figure 4:** Forecasting results of the time series with a linear trend by a) linguistic method; b) linguo-correlation method

As Figure 4 shows, both methods allow to predict dynamical time series, which have a seasonality and linear trend, with high accuracy. In this case, the correlation criterion of similarity showed a slightly better result. It should be noted that studies have not shown a stable advantage of using one or another criterion.

A dynamical time series with a non-linear trend is the most common variant. Based on the available experience, it can be assumed that the use of linear regression to the sequences $X_N$ and $X_K$, and not to the entire time series $X$ in general, will allow us to track not only the linear "global" trends of the series $X$.

**Figure 5** shows the forecasting result of a non-stationary time series of natural numbers with a complex nonlinear trend ($K$=100, $P$=20, $M$ =1, $S$=30) by the linguistic method using linear regression. As expected, the prediction result is quite accurate.

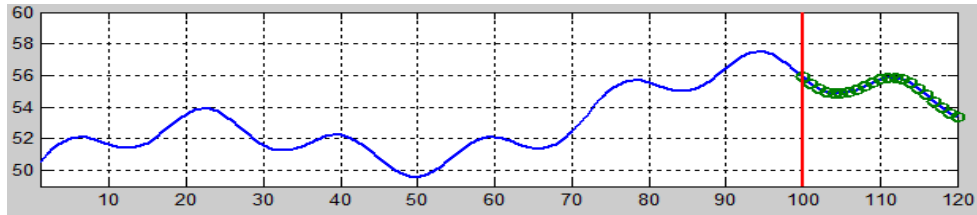

**Figure 5:** Forecasting result of a non-stationary time series with a nonlinear trend by the linguistic method using linear regression

**Figure 6** shows the forecasting result of a similar time series (as the case presented in Figure 5 by the proposed linguistic method, however, this result is noticeably worse. The reason for the deterioration of the forecast accuracy in this case is the lack in the past "phrase" $X_N$ identical $X_K$, since the nonlinear trend is not fully "manifested" within the time series $X$. However, the short-term forecast (at $P \leq 10$) in the case shown in Fig. 6, is also quite accurate.

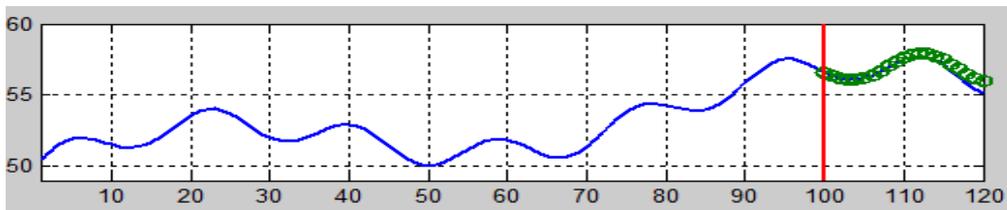

**Figure 6:** Another example of the forecasting result by the mentioned above method

For comparison with other common forecasting methods, **Figure 7** shows the forecasting results of the time series, shown in Figure 6, by the Holt's model. As shown in **Figure 7**, the Holt's model can only be used to obtain a short-term forecast – $P \leq 5$ (Fig. 7 a) or to form "global" trend for rough long-term forecasts (Fig. 7 b).

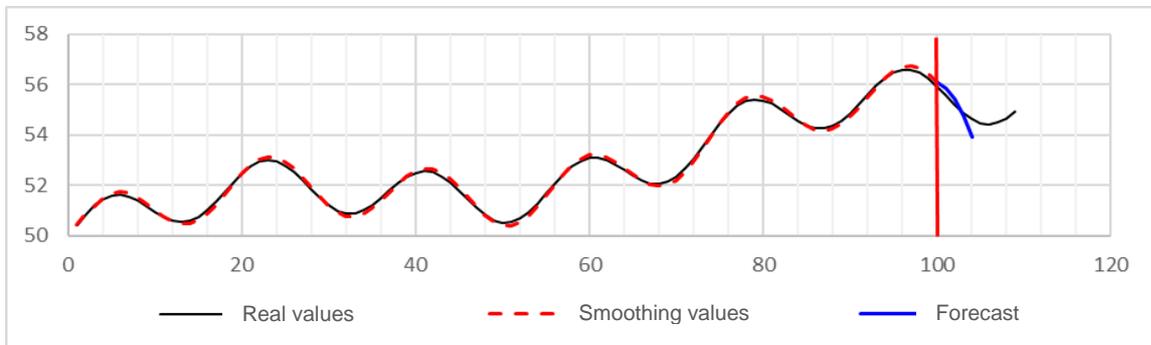

a)

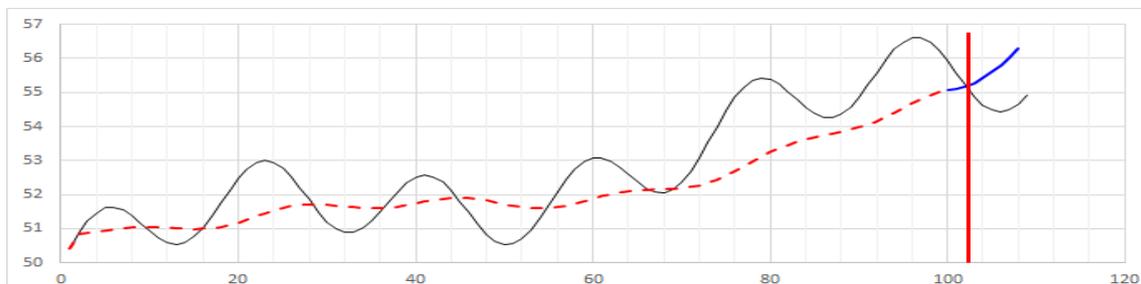

b)

**Figure 7:** Forecasting result of a non-stationary time series with a nonlinear trend by the Holt's model

It remains to be seen whether the proposed linear trends search in the $X_N$ and $X_K$ series of samples by the linear regression method will not lead to errors in the forecast, when in fact there is no global trend. If such impact is significant, an additional step in the analysis of the global time series trend should be added to the linguistic and linguo-correlation methods before forecasting. This, in turn, can offset the advantages of the proposed methods and reduce their versatility. However, performed studies have refuted this warning.
Let's show it by an example. Consider a case similar to that shown in Figure 2, where the time series does not contain a global trend. As can be seen from

**Figure 8**, finding and taking into account the trends of samples $X_N$ and $X_K$ by constructing a linear regression, when there is not a global trend, did not lead to a significant deterioration of the forecast.

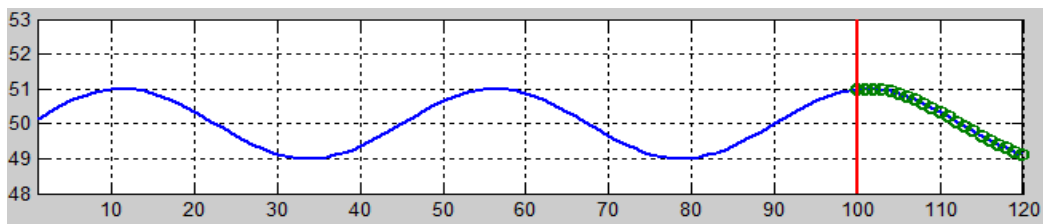

**Figure 8:** Forecasting results of a time series without trend by the linguo-correlation method using linear regression

Figure 9 shows the result of practical application of the linguo-correlation forecasting method of the time series to calculate predicted values of the number of publications in Russian-language media with a mention of the President of United States on the base of the time series (in increments of one day) obtained by the content monitoring system [1,9].

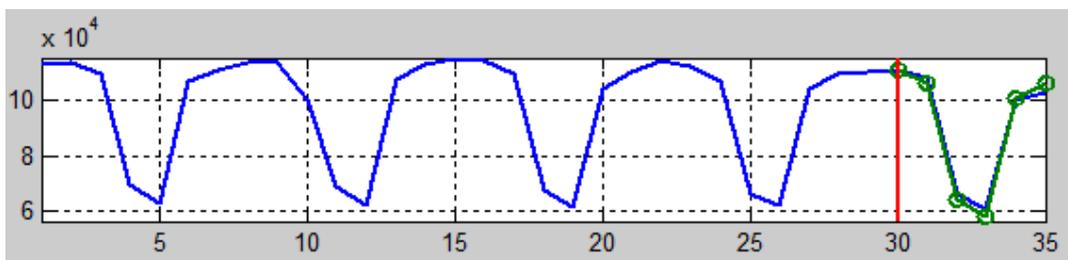

Figure 9: Forecasting result of the publications number in Russian language media with a mention of President of United States by proposed linguistic method

As can be seen from
Figure 9, the predicted values of the series fairly accurately reflect the dynamics of changes in the publications frequency on a given topic and generally correspond to the actual control values.

Thus, the proposed methods of linguistic (without global trend) and linguo-correlation (with trend search) forecasting of the time series have demonstrated their effectiveness and applicability for predicting various dynamic non-stationary time series.

## 4. Conclusions

Proposed methods of linguistic and linguo-correlation forecasting of dynamic time series based on the linguistic approach, the use of the N-gram method and linear regression to analyze local trends of time series belong to the class of autoregressive methods that are easily automated because they have only two input parameters that do not require precise adjustment.

This paper presents examples of application of these methods, shows their features.

Proposed methods provide forecasting for short-term and medium-term periods of dynamic non-stationary time series with both linear and nonlinear trends. Also with a high level of automation they allow to make predictions of time series, which are characterized by cyclicality, in particular, the series of dynamics of publications in content monitoring systems, without additional analysis of their seasonality.

A significant advantages of this approach are the lack of requirements for the time series stationarity and a small number of settings.

Further research may focus on the study of various criteria for the similarity of time series fragments, the use of nonlinear similarity criteria (eg, wave), the search for ways to automatically determine the rational step of the time series quantization.